\documentclass[12pt]{article}
\usepackage[english]{babel}
\usepackage{amssymb}
\usepackage{amsmath}
\usepackage{graphicx}
\usepackage{bbm}
\usepackage[right]{eurosym}

\def\vs{\vspace}
\def\noi{\noindent}

\def\IN{\mathbb N}

\def\IR{\mathbb R}

\def\IQ{\mathbb Q}

\def\IS{\mathbb S}

\def\an{\mathrm{an}}
\def\exp{\mathrm{exp}}

\def\ma{\mathcal}

\def\arg{\mathrm{arg}}

\begin{document}
	
	\begin{center}
		{\bf \Large Growth of 
			Log-Analytic Functions}
	\end{center}
	
	%\vspace{0.3cm}
	\centerline{Tobias Kaiser}
	
	\vspace{0.3cm}\noi 
	\begin{center}
		\begin{minipage}[t]{8cm}\scriptsize{{\bf Abstract.}
				We show that unary log-analytic functions are polynomially bounded.
				In the higher dimensional case globally a log-analytic function can have exponential growth. We show that a log-analytic function is polynomially bounded on a definable set which contains the germ of every ray at infinity.}
		\end{minipage}
	\end{center}
	
	\normalsize
	\section*{Introduction}
	
	Log-analytic functions have been defined by Lion and Rolin in their seminal paper [5]. 
	They are iterated compositions from either side of globally subanalytic functions (see [2]) and the global logarithm. In [3] it was shown that, from the point of view of differentiability, log-analytic functions behave similarly to globally subanalytic functions. We have strong quasianalyticity, and Tamm's theorem holds. But with respect to growth properties, log-analytic functions behave in a different way compared to globally subanalytic functions.  Globally subanalytic functions are polynomially bounded. This holds also for log-analytic funtions of one variable. But in higher dimension, surprisingly, the situation changes. 
	Although the global exponential function is not involved in the definition of log-analytic functions, a log-analytic function in at least two variables can have exponential growth. We construct an example where the function is not polynomially bounded on every dense definable set. But polynomially boundedness holds on a definable set which is `thick' at infinity: We show that a log-analytic function is polynomially bounded on a definable set which contains the germ of every ray at infinity.

	\section*{Notations}

	\vs{0.1cm}
	By $\IN=\{1,2,\ldots\}$ we denote the set of natural numbers, and by $\IN_0=\{0,1,2,\ldots\}$ the set of nonnegative integers.
	
	For $t\in \IR$ we set $\IR_{>t}:=\{x\in\IR\mid x>t\}$ and $\IR_{\geq t}:=\{x\in \IR\mid x\geq t\}$. 
	Denoting by $|\;|$ the euclidean norm on $\IR^n$, we set $\IS^{n-1}:=\{x\in \IR^n\mid |x|=1\}$.
	Given a subset $A$ of $\IR^n$, we denote by $\overline{A}$ its closure.
	
	By $\pi:\IR^n\times\IR\to \IR^n, (x,y)\mapsto x,$ we denote the projection on all but the last coordinate.
	For a subset $A$ of $\IR^n\times \IR$ and $x\in \IR^n$, we set $A_x:=\{y\in \IR\mid (x,y)\in A\}$.
	
	By $\exp_k$, respectively $\log_k$, we denote the $k$-times iterated of the exponential function, respectively the logarithm.

	\vs{0.5cm}
	\hrule
	
	\vs{0.4cm}
	{\footnotesize{\itshape 2010 Mathematics Subject Classification:} 03C64, 14P15, 26A09, 26A12, 32B20}
	\newline
	{\footnotesize{\itshape Keywords and phrases:} log-analytic functions, polynomially bounded, exponential growth}
	
	\newpage

	\section*{The Results}
	
	We assume basic knowledge of o-minimality (see for example van den Dries [1] and van den Dries and Miller [2]). By {\it definable} we mean definable (with parameters) in the o-minimal structure $\IR_{\an,\exp}$ (see [2] for this structure).

	\subsection*{Setting and Preliminaries}
	
	We recall the precise definition of a log-analytic function (see Lion and Rolin [5]) and state consequences of preparation results on special sets (compare with [3]).
	
	\vs{0.5cm}
	{\bf 1. Definition} 
	
	\vs{0.1cm}
	Let $X \subset \IR^n$ be definable and let $f:X \to \IR$ be a function. 
	\begin{itemize}
		\item[(a)] Let $k \in \IN_0$. By induction on $k$ we define that $f$ is {\bf log-analytic of order at most $k$}.
		
		\vs{0.2cm}
		{\bf Base case:} 
		The function $f$ is log-analytic of order at most $0$ if $f$ is piecewise the restriction of globally subanalytic functions; i.e., there is a finite decomposition $\ma{Y}$ of $X$ into definable sets such that for $Y \in \ma{Y}$ there is a globally subanalytic function $F:\IR^n \to \IR$ such that $f|_Y=F|_Y$.
		
		\vs{0.2cm}
		{\bf Inductive step:}
		The function $f$ is log-analytic of order at most $k$ if the following holds: There is a finite decomposition $\ma{Y}$ of $X$ into definable sets such that for $Y \in \mathcal{Y}$ there are $p,q \in \IN_0$, a globally subanalytic function $F:\IR^{p+q} \to \IR$ and log-analytic functions $g_1,...,g_p: Y \to \IR, h_1,\ldots,h_q:Y\to \IR_{>0}$ of order at most $k-1$ such that
		$$f|_Y=F\big(g_1,...,g_p,\log(h_1),...,\log(h_q)\big).$$
		\item[(b)] Let $k\in \IN_0$. We call $f$ {\bf log-analytic of order $k$} if $f$ is log-analytic of order at most $k$ but not of order at most $k-1$.
		\item[(c)] We call $f$ {\bf log-analytic} if it is log-analytic of order $k$ for some $k\in \IN_0$.  
	\end{itemize}
	
	\vs{0.2cm}
	{\bf 2. Definition}
	
	\vs{0.1cm}
	We call a definable cell $Y\subset \IR^{n+1}$ {\bf simple at infinity} if for every $x\in \pi(Y)$ we have $Y_x=\IR_{>d_x}$ for some $d_x\in \IR_{\geq 0}$. 
	
	\vs{0.5cm}
	{\bf 3. Remark}
	
	\vs{0.1cm}
	Let $\ma{Y}$ be a definable cell decomposition of $\IR^n\times\IR_{>0}$. Then
	$$\IR^n=\bigcup\{\pi(Y)\mid Y\in\ma{Y}\mbox{ simple at infinity}\}.$$
	
	\vs{0.5cm}
	We set $e_0:=0$ and $e_k:=\exp(e_{k-1})$ for $k\in \IN$. 
	
	\vs{0.5cm}
	{\bf 4. Definition}
	
	\vs{0.1cm}
	Let $k\in \IN_0$. A cell $Y\subset \IR^{n+1}$ which is simple at infinity is called {\bf $k$-simple at infinity} if $\inf Y_x\geq e_k$ for all $x\in \pi(Y)$.
	
	\vs{0.5cm}
	{\bf 5. Proposition}
	
	\vs{0.1cm}
	{\it Let $f:\IR^n\times\IR\to \IR, (x,y)\mapsto f(x,y),$ be log-analytic of order $k$. Then there is a definable cell decomposition $\ma{Y}$ of $\IR^n\times \IR$ such that for every $Y\in \ma{Y}$ which is simple at infinity the cell $Y$ is $k$-simple at infinity and
		$$f|_Y(x,y)=a(x)y^{q_0}\log(y)^{q_1}\cdots \log_k(y)^{q_k}u(x,y)$$
		where 
		\begin{itemize}
			\item[(1)] $a:\pi(Y)\to \IR$ is log-analytic and continuous,
			\item[(2)] $q_0,\ldots,q_k\in \IQ$,
			\item[(3)] $u:Y\to \IR$ is log-analytic and there is $d\in \IR_{>0}$ such that $0\leq u(x,y)\leq d$ for all $(x,y)\in Y$.
	\end{itemize}}
	{\bf Proof:}
	
	\vs{0.1cm}
	This follows from [3, Theorem 2.30] using the substitution $r\mapsto 1/r$.
	\hfill$\blacksquare$

	\subsection*{Statement and Proof of the Results}
	
	{\bf 6. Definition}
	
	\vs{0.1cm}
	Let $n\in \IN$ and let $f:\IR^n\to \IR$ be a function.
	\begin{itemize}
		\item[(a)] If $n=1$ we say that $f$ is {\bf polynomially bounded at infinity} if there are constants $t\in \IR_{>0}$ and $N\in \IN$ such that
		$|f(x)|\leq x^N$ for all $x>t$.
		\item[(b)]  If $n>1$ we say that $f$ is {\bf polynomially bounded at infinity} if there are constants $t\in \IR_{>0}$ and $N\in \IN$ such that
		$|f(x)|\leq |x|^N$ for all $|x|>t$.
	\end{itemize}
	
	\vs{0.3cm}
	Let $f$ be as above and let $A\subset \IR^n$ be unbounded. We say that $f$ is polynomially bounded at infinity on $A$ if $\mathbbm{1}_A f$ is polynomially bounded at infinity (where $\mathbbm{1}_A f$ denotes the characteristic function of $A$).
	
	\vs{0.5cm}
	We handle the unary case first.
	
	\vs{0.5cm}
	{\bf 7. Proposition}
	
	\vs{0.1cm}
	{\it Let $f:\IR\to \IR$ be log-analytic. Then $f$ is polynomially bounded.}
	
	\vs{0.1cm}
	{\bf Proof:}
	
	\vs{0.1cm}
	By Proposition 5 we find $k\in \IN_0$ and $t\geq e_k$ such that $$f(x)=ax^{q_0}\log(x)^{q_1}\cdots\log_k(x)^{q_k} u(x)$$
	on $\IR_{\geq t}$
	where
	\begin{itemize}
		\item[(1)] $a\in \IR$,
		\item[(2)] $q_0,\ldots,q_k\in \IQ$,
		\item[(3)] $u:\IR_{>t}\to \IR$ is log-analytic and there is $d\in \IR_{>0}$ such that $0\leq u(x)\leq d$ for all $x>t$.
	\end{itemize}
	This gives that $f(x)$ behaves asymptotically as $x^{q_0}\log(x)^{q_1}\cdots\log_k(x)^{q_k}$ at $+\infty$ (unless in the trivial case $a=0$). By the growth properties of the logarithm we are done. \hfill$\blacksquare$
	
	\vs{0.5cm}
	{\bf 8. Definition}
	
	\vs{0.1cm}
	A subset $\ma{C}$ of $\IR^n$ is called a {\bf cone} if $x\in \ma{C}$ implies $rx\in \ma{C}$ for all $r\in \IR_{\geq 0}$.
	
	\vs{0.5cm}
	Given a cone $\ma{C}$ with $\ma{C}\supsetneq \{0\}$ we denote by $B(\ma{C}):=\ma{C}\cap \IS^{n-1}$ its {\bf base}. Note that $\ma{C}=\IR_{\geq 0}\cdot B(\ma{C})$.
	
	\vs{0.5cm}
	{\bf 9. Proposition}
	
	\vs{0.1cm}
	{\it Let $n\geq 2$ and let $f:\IR^n\to \IR$ be log-analytic. Then there is a cone $\ma{C}$ with nonempty interior such that $f$ is polynomially bounded at infinity on $\ma{C}$.}
	
	\vs{0.1cm}
	{\bf Proof:}
	
	\vs{0.1cm}
	We consider the polar coordinates
	$\varphi:\IS^{n-1}\times \IR_{\geq 0}\to \IR^n, (v,r)\to rv$.
	Let $g:\IS^{n-1}\times \IR_{\geq 0}\to \IR, (v,r)\mapsto f(\varphi(v,r))$.
	By Remark 3 and Proposition 5  we find $k\in \IN_0$ and an open cell $Y$ that is $k$-simple at infinity such that $g|_Y(x)=a(v)r^{q_0}\log(r)^{q_1}\cdots\log_k(r)^{q_k} u(v,r)$
	where
	\begin{itemize}
		\item[(1)] $a:\pi(Y)\to \IR$ is log-analytic and continuous,
		\item[(2)] $q_0,\ldots,q_k\in \IQ$,
		\item[(3)] $u:Y\to \IR$ is log-analytic and there is $d\in \IR_{>0}$ such that $0\leq u(v,r)\leq d$ for all $(v,r)\in Y$.
	\end{itemize}
	Choose an open ball $B$ in $\pi(Y)$ such that its closure is contained in $\pi(Y)$.
	Then by continuity $a$ is bounded on $B$.
	By the growth properties of the iterated logarithms, we get that $g$ is polynomially bounded on $Y\cap (B\times \IR)$. 
	By the definition of cells, the map $\pi(Y)\to \IR_{\geq 0}, x\mapsto \inf Y_x,$ is continuous. Hence by the conditions imposed on $B$,
	there is $T>e_k$ such that the function $x\mapsto \inf Y_x$ on $B$ is bounded from above by $T$. This implies $B\times \IR_{>T}\subset Y$. Hence, $g$ is polynomially bounded on $B\times \IR_{>T}$. We consider the cone 
	$\ma{C}:=\IR_{\geq 0}\cdot B$ which has nonempty interior. We obtain some $N\in \IN$ such that $|f(x)|\leq |x|^N$ for all $x\in \ma{C}$ with $|x|>T$. By the very definition we obtain that $f$ is polynomially bounded at infinity on $\ma{C}$.
	\hfill$\blacksquare$
	
	\vs{0.5cm}
	In the higher dimensional case, global (polynomially) boundedness may fail simply if the pole locus is not bounded.
	Consider for example the function
	$$f:\IR^2\to \IR, (x,y)\mapsto \left\{\begin{array}{ccc}
		\frac{1}{x-y}&&x\neq y,\\
		&\mbox{if}&\\
		0&&x=y.\end{array}	\right.$$
	Then clearly $\sup_{\sqrt{x^2+y^2}=r}|f(x,y)|=\infty$ for all $r>0$.
	
	But even if one restricts to continuous functions, a log-analytic function may be not be polynomially bounded if $n\geq 2$.
	
	\vs{0.5cm}
	{\bf 10. Proposition}
	
	\vs{0.1cm}
	{\it Let $n\geq 2$. There is a continuous log-analytic function $f:\IR^n\to \IR$ which is not polynomially bounded at infinity.}
	
	\vs{0.1cm}
	{\bf Proof:}
	
	\vs{0.1cm}
	It suffices to deal with the case $n=2$. 
	Consider the function
	$$h:\IR_{>1}\times\IR_{>0}\to \IR, (x,y)\mapsto -y\big((\log y)^2-2\log y+2-x\big).$$
	
	\vs{0.2cm}
	{\bf Claim 1:}
	The following holds:
	\begin{itemize}
		\item[(1)] The function $h$ is log-analytic and continuous.
		\item[(2)] For every $x>1$ there exists $\max_{y>0}h(x,y)\in \IR$.
		\item[(3)] The function $\alpha:\IR_{>1}\to \IR, x\mapsto \max_{y>0}h(x,y),$ is given by $\alpha(x)=2\exp(\sqrt{x})(\sqrt{x}-1)$.
	\end{itemize}
	
	{\bf Proof of Claim 1:}
	
	\vs{0.1cm}
	(1) being clear, we have to show (2) and (3). For $x>1$ we have
	$$\lim_{y\nearrow \infty}h(x,y)=-\infty, \lim_{y\searrow 0} h(x,y)=0$$
	and
	$$\frac{\partial h}{\partial y}(x,y)=-(\log y)^2+x$$
	which vanishes exactly for $y=\exp(\sqrt{x})$ and $y=\exp(-\sqrt{x})$.
	We have
	$$h(x,\exp(\sqrt{x}))=2\exp(\sqrt{x})(\sqrt{x}-1), h(x,\exp(-\sqrt{x}))=-2\exp(\sqrt{x})(\sqrt{x}+1).$$
	This implies that for $x>1$ the function $\IR_{>0}\to \IR, y\mapsto h(x,y),$ attains its maximum at $y=\exp(\sqrt{x})$ with this maximum being given by 
	$$\max_{y>0}h(x,y)=2\exp(\sqrt{x})(\sqrt{x}-1).$$
	This shows (2) and (3).
	\hfill$\blacksquare_{\mathrm{Claim\,1}}$
	
	\vs{0.2cm}
	Let $a\in \IR_{>1}$ be the (uniquely determined) value such that $2\exp(\sqrt{a})(\sqrt{a}-1)=1$. 
	Let	
	$$g:\IR_{\geq 0}\times [0,1]\to\IR,(x,y)\mapsto\left\{
	\begin{array}{ccc}
		\max\big\{h\big(x,\frac{y}{1-y}),1\big\},&&(x,y)\in \IR_{>a}\times \,]0,1[,\\
		&\mbox{if}&\\
		1,&&(x,y)\notin \IR_{>a}\times \,]0,1[.
	\end{array}\right.$$
	
	\vs{0.2cm}
	{\bf Claim 2:}
	The following holds:
	\begin{itemize}
		\item[(1)] The function $g$ is continuous and log-analytic.
		\item[(2)] The function $\beta:\IR_{\geq 0}\to \IR, x\mapsto \max_{0\leq y\leq 1}g(x,y),$ is given by $\beta(x)=1$ for $x\leq a$ and $\beta(x)=\alpha(x)$ for $x>a$.
	\end{itemize}
	
	{\bf Proof of Claim 2:}
	
	\vs{0.1cm}
	For (1) note that for $b>a$ 
	$$\lim_{x\to b, y\nearrow 1}g(x,y)=\lim_{x\to b, y\nearrow \infty}\max\big\{h(x,y),1\big\}=1,$$ 
	$$\lim_{x\to b, y\searrow 0}g(x,y)=\lim_{x\to b, y\searrow 0}\max\big\{h(x,y),1\big\}=1,$$
	that for $0<c<1$
	$$\lim_{x\searrow a, y\to c}g(x,y)=\lim_{x\searrow a, y\to c}\max\big\{h(x,y/(1-y)),1\big\}=1,$$
	and that
	$$\lim_{x\searrow a, y\searrow 0}g(x,y)=\lim_{x\searrow a, y\searrow 0}\max\big\{h(x,y/(1-y)),1\big\}=1,$$
	$$\lim_{x\searrow a, y\nearrow 1}g(x,y)=\lim_{x\searrow a, y\nearrow 1}\max\big\{h(x,y/(1-y)),1\big\}=1.$$
	For (2) note that for $x>a$
	$$\max_{0\leq y\leq 1}g(x,y)=\max_{y>0}h(x,y)=2\exp(\sqrt{x})(\sqrt{x}-1).$$
	\hfill$\blacksquare_{\mathrm{Claim\,2}}$
	
	\vs{0.2cm}
	Let
	$$f:\IR^2\to \IR,
	(x,y)\mapsto \left\{
	\begin{array}{ccc}
		g\Big(x^2+y^2, \frac{1}{2\pi}\arg\big(\frac{(x,y)}{\sqrt{x^2+y^2}}\big)\Big),&&(x,y)\neq (0,0),\\
		&\mbox{if}&\\
		1,&&(x,y)=(0,0),
	\end{array}\right.$$
	where the argument function is given by
	$\arg:\IS^1\to [0,2\pi[$  with $\arg((1,0))=0$ and counterclockwise orientation. 
	Then $f$ is continuous and log-analytic.
	Let $\gamma:\IR_{\geq 0}\to \IR_{\geq 0}, r\mapsto \max_{\sqrt{x^2+y^2}=r}|f(x,y)|$.
	Then $\gamma(r)=\alpha(r^2)$ for all $r\geq 0$. Hence,
	$$\mathrm{max}_{\sqrt{x^2+y^2}=r}|f(x,y)|\geq \exp(r)$$
	for all sufficiently large $r$.
	\hfill$\blacksquare$
	
	\vs{0.5cm}
	The question is how ``big'' we can choose a set where polynomially boundedness at infinity holds. In Proposition 9 we have shown that we can choose a nonempty open cone.
	By the continuity of the counterexample in Proposition 10 we cannot hope for a dense definable set (or equivalently, a definable set with dimension of the complement being smaller than $n$):
	
	\vs{0.5cm}
	{\bf 11. Corollary}
	
	\vs{0.1cm}
	{\it Let $n\geq 2$. There is log-analytic function $f:\IR^n\to \IR$ such that $f$ is not polynomially bounded on every dense definable subset.}
	
	\vs{0.5cm}
	{\bf 12. Remark}
	
	\vs{0.1cm}
	Note that the above counterexample is globally given by composition of globally subanalytic functions and the logarithm, not only piecewise.
	
	\vs{0.5cm}
	To formulate an optimal result we need to introduce some setting to speak about the ultimate size of a set at $\infty$.
	The first definition mimics the tangential cone at finite points (see for example Kurdyka and Raby [4]).
	
	We fix an unbounded definable subset $A$ of $\IR^n$. We let $\dim_\infty A$ be $\dim (A\cap \{x\in \IR^n\mid |x|>r\})$ for sufficiently large $r$ (note that this stabilizes) and call it the {\bf dimension} of $A$ {\bf at infinity}.
	
	\vs{0.5cm}
	{\bf 13. Definition} 
	
	\begin{itemize}
		\item[(a)]
		We let 
		$B(A,\infty)$ be the set of all $v\in \ma{S}^{n-1}$ such that for every $r,\varepsilon>0$ there is $x\in A$ with $|x|>r$ and $\big\vert x/|x|-v\big\vert<\varepsilon$.
		We call $\ma{C}(A,\infty):=\IR_{\geq 0}\cdot B(A,\infty)$ the {\bf tangent cone} of $A$ {\bf at infinity}. 
		\item[(b)]
		We let 
		$B^\mathrm{str}(A,\infty)$ be the set of all $v\in \ma{S}^{n-1}$ such that there is some $t\in \IR_{\geq 0}$ with $\IR_{\geq t}\cdot v\subset A$.
		We call $\ma{C}^\mathrm{str}(A,\infty):=\IR_{\geq 0}\cdot B^\mathrm{str}(A,\infty)$ the {\bf strong tangent cone} of $A$ {\bf at infinity}.
	\end{itemize}
	
	\vs{0.5cm}
	{\bf 14. Remark}
	
	\begin{itemize}
		\item[(1)]
		We have $\ma{C}^\mathrm{str}(A,\infty)\subset\ma{C}(A,\infty)$. 
		\item[(2)]
		The tangent cone $\ma{C}(A,\infty)$ of $A$ at infinity is closed and definable with $\dim \ma{C}(A,\infty)\leq \dim_\infty A$. 
		\item[(3)]
		The strong tangent cone $\ma{C}^\mathrm{str}(A,\infty)$ of $A$ at infinity is definable with $\dim \ma{C}^\mathrm{str}(A,\infty)\leq \dim_\infty A$.
		\item[(4)]
		For $r>0$ let $B(A,r):=\{x/r\mid x\in A\mbox{ and }|x|=r\}$. Then $B(A,\infty)$ is the Hausdorff limit of the family $\big(\overline{B(A,r)}\big)_{r\in \IR_{>0}}$ (compare with Lion and Speissegger [6]) and $$B^\mathrm{str}(A,\infty)=\limsup_{r>0}B(A,r)=\bigcup_{r>0}\bigcap_{s>r}B(A,s).$$
	\end{itemize}
	
	\vs{0.5cm}
	The next concept will carry more information.
	A {\bf (closed) ray} $\ma{R}$ in $\IR^n$ is of the form $\ma{R}=a+\IR_{\geq 0}\cdot v$ where $a\in \IR^n$ and $v\in \IS^{n-1}$.
	We parameterize the set $\mathfrak{R}$ of all rays 
	by the bijection
	$\IR^n\times \IS^{n-1}\to \mathfrak{R}, (a,v)\mapsto a+\IR_{\geq 0}\cdot v$.
	For limit considerations, it is natural to identify two rays $\ma{R}_1$ and $\ma{R}_2$ if $\ma{R}_1\subset \ma{R}_2$ or $\ma{R}_2\subset \ma{R}_1$. This is an equivalence relation $\sim$ on $\mathfrak{R}$. A canonical representative of the equivalence class of a ray $\ma{R}=a+\IR_{\geq 0}\cdot v$ is given by $o+\IR_{\geq 0}\cdot v$ where $o\in a+\IR\cdot v$ with
	$o\perp v$ (or, equivalently, $o$ realizes the distance of the line $a+\IR\cdot  v$ to the origin). A ray of this form is called a {\bf standardized ray}.
	We identify the set $\mathfrak{R}/\sim$ with the set of the standardized rays
	and parameterize it by the bijection
	$\mathfrak{S}:=\{(o,v)\in \IR^n\times \IS^{n-1}\mid o\perp v\}\to \mathfrak{R}/\sim,
	(o,v)\mapsto o+\IR_{\geq 0}\cdot v$.

	\vs{0.5cm}
	{\bf 15. Definition} 
	
	\begin{itemize}
		\item[(a)]
		We denote by $\ma{RC}(A,\infty)$ the union of all standardized rays $\ma{R}=o+\IR_{\geq 0}\cdot v$ 
		such that for every $r,\varepsilon>0$, there are $x\in A$ and $y\in \ma{R}$ with $|x|=|y|>r$ and $|x-y|<\varepsilon$, and call it the 
		the {\bf tangent ray cone} of $A$ {\bf at infinity}. 
		\item[(b)]
		We denote by $\ma{RC}^\mathrm{str}(A,\infty)$ the union of all standardized rays $\ma{R}=o+\IR_{\geq 0}\cdot v$ such that $o+\IR_{\geq t}\cdot v\subset A$ for some $t\in\IR_{\geq 0}$, and call it the 
		{\bf strong tangent ray cone} of $A$ {\bf at infinity}. 
	\end{itemize}
	
	\vs{0.5cm}
	{\bf 16. Remark}
	
	\begin{itemize}
		\item[(1)]
		We have $\ma{RC}^\mathrm{str}(A,\infty)\subset \ma{RC}(A,\infty)$.
		\item[(2)]
		The tangent ray cone $\ma{RC}_{A,\infty}$ of $A$ at infinity is closed and definable with $\dim \ma{RC}_{A,\infty}\leq \dim_\infty A$.
		\item[(3)] 
		The strong tangent ray cone $\ma{RC}^\mathrm{str}_{A,\infty}$ of $A$ at infinity is definable with $\dim \ma{RC}^\mathrm{str}_{A,\infty}$
		$\leq \dim_\infty A$.
		\item[(4)]
		We have $\ma{C}(A,\infty)\subset \ma{RC}(A,\infty)$. In fact, the following stronger statement holds: A standardized ray $o+\IR_{\geq 0}\cdot v$ is contained in $\ma{RC}(A,\infty)$ if and only if $\IR_{\geq 0}\cdot v$ is contained in $\ma{C}(A,\infty)$.
		\item[(5)] 
		We have $\ma{C}^\mathrm{str}(A,\infty)\subset \ma{RC}^\mathrm{str}(A,\infty)$.
	\end{itemize}
	
	\vs{0.2cm}
	{\bf 17. Example}
	
	\vs{0,1cm}
	Consider the half-strip
	$$S:=\{(x,y)\in \IR^2\mid x>0, 0<y<1\}.$$
	We have
	$$\ma{C}(S,\infty)=\IR_{\geq 0}\cdot (1,0), \ma{C}^\mathrm{str}(S,\infty)=\emptyset$$ and 
	$$\ma{RC}(S,\infty)=\big\{(0,t)+\IR_{\geq 0}\cdot (1,0)\,\big\vert\; t\in\IR\big\},$$
	$$\ma{RC}^\mathrm{str}(S,\infty)=\big\{(0,t)+\IR_{\geq 0}\cdot (1,0)\;\big\vert\; 0<t< 1\big\}.$$
	
	\vs{0.5cm}
	{\bf 18. Definition}
	
	\begin{itemize}
		\item[(a)]
		We call $A$ {\bf spherically dense at infinity} if $\ma{C}(A,\infty)=\IR^n$. 
		We call $A$ {\bf strongly spherically dense at infinity} if $\ma{C}^\mathrm{str}(A,\infty)=\IR^n$.
		\item[(b)] 
		We call $A$ {\bf ray dense at infinity} if $\ma{RC}(A,\infty)$ contains every standardized ray.
		We call $A$ {\bf strongly ray dense at infinity} if $\ma{RC}^\mathrm{str}(A,\infty)$ contains every standardized ray.
	\end{itemize}
	
	\vs{0.1cm}
	{\bf 19. Remark}
	
	\begin{itemize}
		\item[(1)] $A$ is spherically dense at infinity if and only if $A$ is ray dense at infinity.
		\item[(2)] 
		If $A$ is strongly ray dense at infinity then $A$ is strongly spherically dense at infinity. The converse does in general not hold.
	\end{itemize}
	{\bf Proof:}
	
	\vs{0.1cm}
	(1): The direction from right to the left being clear by definition, we show the direction from left to the right.
	Let $o+\IR_{\geq 0}\cdot v\in \mathfrak{R}/\sim$ where $(o,v)\in \mathfrak{S}$. Then $\IR_{\geq 0}\cdot v\in \ma{C}(A,\infty)$ since $A$ is spherically dense at infinity.
	By the definition of the tangent ray cone, we obtain that $o+\IR_{\geq 0}\cdot v\subset \ma{RC}(A,\infty)$.
	
	\vs{0.2cm}
	(2): The first statement ist clear. For the second one, consider the complement of the above half-strip. 
	\hfill$\blacksquare$
	
	\vs{0.5cm}
	Hence the notion of ray density at infinity does not give anything new. We have included it for completeness and symmetry.
	
	\vs{0.5cm}
	Here is now the final optimal result.
	
	\vs{0.5cm}
	{\bf 20. Theorem}
	
	\vs{0.1cm}
	{\it Let $n\geq 2$ and let $f:\IR^n\to \IR$ be log-analytic. Then there is a definable subset $\ma{U}$ of $\IR^n$ which is strongly ray dense at infinity such that $f$ is polynomially bounded at infinity on $\ma{U}$.}
	
	\vs{0.1cm}
	{\bf Proof:}
	
	\vs{0.1cm}
	Consider the semialgebraic map
	$\Phi:\mathfrak{S}\times \IR_{\geq 0}\to \IR^n, (o,v,r)\mapsto o+rv,$
	and the log-analytic function
	$F:=f\circ \Phi:\mathfrak{S}\times \IR_{\geq 0}\to \IR$. Let $F$ be log-analytic of order $k\in\IN_0$.
	By Proposition 5 we find a definable cell decomposition $\ma{Y}$ of $\mathfrak{S}\times \IR_{\geq 0}$ such that for every $Y\in \ma{Y}$ which is simple at infinity the cell $Y$ is $k$-simple at infinity such that
	$$F|_Y(o,v,r)=a(o,y)r^{q_0}\log(r)^{q_1}\cdots \log_k(r)^{q_k}u(o,v,r)$$
	where 
	\begin{itemize}
		\item[(1)] $a:\pi(Y)\to \IR$ is log-analytic and continuous,
		\item[(2)] $q_0,\ldots,q_k\in \IQ$,
		\item[(3)] $u:Y\to \IR$ is log-analytic, and there is $d=d_Y\in \IR_{>0}$ such that $0\leq u(o,v,r)\leq d$ for all $(o,v,r)\in Y$.
	\end{itemize}
	We fix $Y\in \ma{Y}$ simple at infinity. Let $Z:=\pi(Y)$ and $\delta:Z\to \IR_{\geq 0}, (o,v)\mapsto \inf Y_{(o,v)}$. We set $\mathrm{fr}Z_\mathfrak{S}:=(\overline{Z}\setminus Z)\cap \mathfrak{S}$. By passing to a finer cell decomposition of $\mathfrak{S}$, we may assume that $\mathrm{fr}Z_\mathfrak{S}\neq \emptyset$. 
	For $s\in \IR_{\geq 0}$ let 
	$$Z(s):=\big\{(o,v)\in Z\,\big\vert\; |(o,v)|\leq s, \mathrm{dist}((o,v),\mathrm{fr}_\mathfrak{S} Z)\geq s\big\}.$$
	Then $Z(s)$ is compact for every $s\geq 0$. We set
	$$\Delta:\IR_{\geq 0}\to \IR_{\geq 0}, s\mapsto \max\big\{|a(o,v)|\;\big\vert\; (o,v)\in Z(s)\big\}.$$
	Note that this is well-defined since $a$ is continuous. Note that here by convention $\max \emptyset =0$. 
	The function $\Delta$ is increasing and definable. Hence, by van den Dries and Miller [2, 5.5], it is bounded by an iterated exponential $\exp_l$ for some $l\in \IN_0$.
	Choose $N=N_Y\in \IN$ with $N>|q_0|+\ldots+|q_n|$.
	We set
	$$\ma{W}_Y:=\big\{(o,v,r)\in \mathfrak{S}\times \IR_{>0}\;\vert\;(o,v)\in Z(\log_l(r)), r>\max\{e_l,\delta(o,v)\} \big\}.$$
	For $(o,v,r)\in \ma{W}_Y$ we have
	$$|F(o,v,r)|=|a(o,v)|r^{q_0}\log(r)^{q_1}\cdots\log_k(r)^{q_k}u(o,v,r)
	\leq d_Yrr^{N_Y}.$$
	We set $\ma{V}_Y:=\Phi(\ma{W}_Y)$. We obtain that
	$|f(x)|\leq d_Y |x|^{N_Y+1}$ on $\ma{V}_Y$. 
	
	\newpage
	Let $\ma{U}$ be the union of all $\ma{V}_Y$ with $Y\in \ma{Y}$ simple at infinity. Then $\ma{U}$ is definable.
	We show that this $\ma{U}$ does the job. 
	Let $\ma{R}=o+\IR_{\geq 0}\cdot v$ be a standardized ray and let $r>0$.
	By Remark 3 we find $Y\in \ma{Y}$ that is simple at infinity such that $(o,v)\in Z$. Note that we use the above notations. 
	There is $s\in \IR_{>0}$ such that $(o,v)\in Z(s)$. By the definition of $\ma{W}_Y$ we find $t>0$ such that
	$\{(o,v)\}\times\IR_{\geq t}\subset \ma{W}_Y$. This gives $o+\IR_{\geq t}\cdot v\subset \ma{V}_Y\subset \ma{U}$. So $\ma{U}$ is strongly ray dense.
	Let 
	$$d_\ma{U}:=\max\{d_y\mid Y\in \ma{Y}\mbox{ simple at infinity}\}$$
	and
	$$N_\ma{U}:=\max\{N_y\mid Y\in \ma{Y}\mbox{ simple at infinity}\}.$$
	Then $|f(x)|\leq d_\ma{U}|x|^{N_\ma{U}+1}$ for all $x\in \ma{U}$. Hence, $f$ is polynomially bounded on $\ma{U}$.
	\hfill$\blacksquare$. 
	
	\vs{0.5cm}
	{\bf 21. Concluding Remarks}
	
	\vs{0.1cm}
	In Corollary 11 we have found for $n\geq 2$ a  log-analytic function $f:\IR^n\to \IR$ and a definable open und unbounded set $W$ such that $r\mapsto \inf_{x\in W, |x|=r}|f(x)|$ is of exponential growth. 
	By Proposition 7 the set $W$ cannot contain the image of an unbounded log-analytic curve.
	By the same methods as in the proof of Theorem 20, we can find an open and definable set $\ma{U}$ such that $f$ is polynomially bounded at infinity on $\ma{U}$, and $\ma{U}$ contains the germ of every unbounded log-analytic curves up to a certain complexity (where the complexity is the complexity of terms in the language $\ma{L}_\an(^{-1},(\sqrt[n]{...})_{n=2,3,...},\log)$, compare with [3, Remark 1.2]). 
	An open question is whether we can find such $\ma{U}$ that contains the germ of every unbounded log-analytic curve.

	\vs{1cm}
	\noi \footnotesize{\centerline{\bf References}
		\begin{itemize}
			\item[(1)] 
			L. van den Dries: Tame Topology and O-minimal Structures. {\it London Math. Soc. Lecture Notes Series} {\bf 248}, Cambridge University Press, 1998.
			\item[(2)]
			L. van den Dries and C. Miller:
			Geometric categories and o-minimal structures.
			{\it Duke Math. J.} {\bf 84} (1996), no. 2, 497-540.
			\item[(3)]
			T. Kaiser and Andre Opris: Differentiability Properties of Log-Analytic Functions.
			{\it Rocky Mountain Journal of Mathematics} {\bf 52} (2022) no. 4, 1423-1443.
			\item[(4)] 
			K. Kurdyka, G. Raby:  Densit\'e des ensembles sous-analytiques.
			{\it Ann. Inst. Fourier} {\bf 39} (1989), no. 3, 753-771. 
			\item[(5)]
			J.-M. Lion, J.-P. Rolin:
			Th\'{e}or\`{e}me de pr\'{e}paration pour les fonctions logarithmico-exponentielles.
			{\it Ann. Inst. Fourier} {\bf 47} (1997), no. 3, 859-884.
			\item[(6)] J.-M. Lion, P. Speissegger:
			A geometric proof of the definability of Hausdorff limits. 
			{\it Sel. Math., New Ser.} {\bf 10} (2004), no. 3, 377-390. 
			
	\end{itemize}}
	
	\vs{0.5cm}
	Tobias Kaiser,
	University of Passau,
	Faculty of Computer Science and Mathematics\\
	tobias.kaiser@uni-passau.de,
	D-94030 Germany

\end{document}